\documentclass{amsart}
\usepackage{amssymb}
\usepackage{amsfonts}
\usepackage{amsmath}
\usepackage{lastpage}
\usepackage[legalpaper,bookmarks=true,colorlinks=true,linkcolor=blue,citecolor=blue]{hyperref}
\usepackage{graphicx}%
\usepackage{fancyhdr}
\usepackage{color}
\usepackage[mathlines]{lineno}
\usepackage{lscape}
\usepackage{epsfig}
\usepackage{natbib}
\usepackage{geometry}
\usepackage{tgbonum}
\fontfamily{qcr}\selectfont

\newtheorem{theorem}{Theorem}
\theoremstyle{plain}

\newtheorem{fact}{Fact}

\newtheorem{definition}{Definition}

\newtheorem{proposition}{Proposition}

\numberwithin{equation}{section}

\newcommand{\Bin}{\bigskip \noindent}

\newcommand{\Ni}{\noindent}

\begin{document}
	\Large

\title[A modern approach to the moment problem on $\mathbb{R}$]{A modern approach to the moment problem on $\mathbb{R}$}
\author{Moussoda Tour\'e $^{\dag}$} 
\author{Gane Samb Lo $^{\dag\dag}$}
\author{Aladji Babacar Niang $^{\dag\dag\dag}$}

\begin{abstract} The moment problem is an important problem in Functional Analysis and in Probability measure. It goes back to Stieltjes, around 1890. There is still an important ongoing interest in the recent literature. But, up today, the main theoretical resource (Shohat and Tamarkin, 1934) does not have the modern exposure it deserves, especially in the current development of measure theory of integration. Besides, the multivariate version is far less exploited. In this paper, a full exposure of such a theory is presented, using the latest knowledge of measure theory and functional analysis. As a result, the basis of future development is layed out and the accessibility of the theory by modern graduate students and researches is guaranteed.\\

\noindent Moussoda Tour\'e $^{\dag}$\\
LERSTAD, Gaston Berger University, Saint-Louis, S\'en\'egal.\\
moussodatoure2014@gmail.com \\

\noindent $^{\dag}$ Gane Samb Lo $^{\dag\dag}$.\\
LERSTAD, Gaston Berger University, Saint-Louis, S\'en\'egal (main affiliation).\newline
LSTA, Pierre and Marie Curie University, Paris VI, France.\newline
AUST - African University of Sciences and Technology, Abuja, Nigeria\\
gane-samb.lo@edu.ugb.sn, gslo@aust.edu.ng, ganesamblo@ganesamblo.net\\
Permanent address : 1178 Evanston Dr NW T3P 0J9,Calgary, Alberta, Canada.\\

\noindent Aladji Babacar Niang $^{\dag\dag\dag}$\\
LERSTAD, Gaston Berger University, Saint-Louis, S\'en\'egal.\\
Email: aladjibacar93@gmail.com\\

\end{abstract}

\maketitle

\section{Introduction}

\bigskip\noindent The problem of moment  is an interesting topic in Functional Analysis, especially in measure theory. It has important applications in probability theory.\\

\noindent Although there is a significant number of research works in probability theory on this problem (see \cite{gutt}, \cite{billinsgleymp},\cite{loeve} and references therein, etc.), the most  important source of that question, when treated in its generality, is \cite{shohat}. Up to our knowledge, we did not see another full set up of that theory beyond that main reference.\\

\noindent We already pointed out that this problem is used in Probability Theory, but the following special form : given a probability law $\mathbb{P}_X$on $\mathbb{R}$ having moments of all orders $(m_n)_{n\geq 0}$, does the sequence $(m_n)_{n\geq 1}$ uniquely determine the probability law $\mathbb{P}_X$? This is a consequence of the moment problem, which goes back to Stieltjes (see \cite{shohat} for references on all particular form of that problem) formulated as follows :\\

\Bin (\textbf{Stieltjes's problem}) [Around 1890, see \cite{shohat} and references therein]. Given $(m_n)_{n\geq 1} \subset \mathbb{R}_{+}$, does it exists finite measure $\rho$ supported by $\mathcal{V}=\mathbb{R}_{+}$ such that

\begin{equation}
\forall n\geq 0, \ m_n=\int_{\mathcal{V}} x^n \ d\rho(x). \label{mp_01}
\end{equation} 

\noindent Of course, $m_0=\rho(\mathbb{R})=\rho(\mathcal{V})$. The term $m_0\neq 0$ is the bound of $\rho$. Later, the same problem is set for a general sequence of real numbers and for $\mathcal{V}=\mathbb{R}$ or  $\mathcal{V}=[0,1]$ and is named after Hamburger and Hausdorff respectively.\\

\Bin The general solution of the problem, when the support $\mathcal{V}$ is bounded by a closed set $S_0$, is given in \cite{shohat}. From there, we face two major concerns about the exposition of the general theory.\\

\bigskip\noindent First, the paper of \cite{shohat}, in our view, is written with the Stieltjes integrals and is based on the rudimentary tools of measure theory and weak convergence of that time of 1943. During the preparation of a master degree dissertation of the second author, we find out that a lot of arguments used by \cite{shohat} may be replaced by arguments that are common now and more appropriate. Essentially, the authors used the Stieltjes integration, the notion of \textit{substancially} continuity or of \textit{substancially} convergences, extension theorems, etc., all those tools seeming to be obsolete now.\\

\noindent By using the modern Lebesgue-Stieltjes integration, the modern theory of weak convergence, the extension theorems of measures on semi-algebras or on algebras, the Caratheodory theorem instead for example, in one word, measure theory arguments, we think that this master-piece paper on the topic can be rendered into a far more readable text for mathematicians of our modern days.\\

\noindent Secondly, the proofs of \cite{shohat} are directly given on $\mathbb{R}^d$, $d\geq 1$. By comparison, classical graduate textbooks in probability refer to the moment problem in one dimension and common readers are used to a multivariate approach of the problem of moments.\\

\noindent Based on the importance of the question and its connections to the characterizations of the weak convergence through the convergence of the moments (it they all exist), we wish to produce a general introduction to the question and entirely expose it at the light of the modern theorem under the following organization :\\

\noindent (1) Treating entirely the dimensional stage with the full details and address the weak convergence through the convergence of the moments (as in 
\cite{billinsgleymp} and \cite{loeve}).\\


\noindent (2) By exposing the ideas of \cite{shohat}, our contribution is two-fold :\\

\noindent (2a) We provide relevant complements and variety of modern arguments that will make the text readable just after a course of Measure Theory and Probability Theory. We intend to formulate the main theorem in \cite{shohat} in the frame of measure theory with the help of some well-known criteria. But, at least, we include needed the  mathematical background. At the end, we hope that a graduate student will be able to read it more comfortably.\\ 

\noindent (2b) In the proofs themselves, we bring more clarity on the linear spaces on which the linear mapping is constructed (see Step 1 in page \pageref{step1}). In the original paper, the roles of $r$ is ambiguous. Actually, the right space should be the class of functions bounded by finite linear combinations of functions 
$u\mapsto A(u_1^{2r_i}+\cdots+u_d^{2r_d})+B$, $A\geq 0$ and $B\geq 0$ (in dimension $d\geq 1$) with non-negative coefficients.\\

\noindent (2c) All along the proof, the right modern tool is used, in particular the Fatou-Lebesgue theorem and the construction of the Lebesgue definition for measurable function of constant sign.\\

\bigskip\noindent Let us organize the paper as follows.\\

\noindent In the next section \ref{mp_maths}, we state the tools we are going to use on modern theory of distribution functions, Lebesgue-Stieltjes integration, limit theorems, etc.\\

\noindent In the Section \ref{mp_r_proba}, we deal with the moment problem within Probability Theory on $\mathbb{R}$ and link it to weak convergence, following mainly \cite{billinsgleymp}.\\

\noindent In Section \ref{mp_r_shohat}, we expose the full proof of \cite{shohat} on $\mathbb{R}$.\\




\newpage
\section{Mathematical background} \label{mp_maths}

\noindent \textbf{A - Distribution functions on $\mathbb{R}^d$, $d\geq 1$}.\\

\noindent The properties we summarize in this Part can be found in major sources as \cite{loeve}, \cite{billingsley}, etc. or in \cite{ips-mestuto-ang} (Chapter 11, 
page 664) for the links between distribution function and Lebesgue-Stieltjes measures and in \cite{ips-wcia-ang} for $F$-continuous intervals.\\
 
\noindent \textbf{A1- Recalls of definitions}. Let us introduce the following internal operation on  $\mathbb{R}^d$:

\begin{equation}
(x,y)\ast (X,Y)=(x_1X_1, x_2X_2, ..., y_dY_d). \label{prodRd}
\end{equation}

\bigskip\noindent Let us consider a real-valued function $F$, defined as follows:

\begin{equation*}
\begin{array}{ccc}
\mathbb{R}^{d} & \mapsto  & \mathbb{R} \\
t & \hookrightarrow  & F(t).
\end{array}
\end{equation*}

\noindent For any interval of $\mathbb{R}^d$ of the form

$$
]a,b]=\prod_{i=1}^{d} ]a_i,b_i]
$$

\noindent for $a=(a_{1},...,a_{d})\leq b=(b_{1},...,b_{d})$, in the sense that $a_i\leq b_i$ for all $i\in \{1,\cdots,d\}$, we define  its $F$-volume by

\begin{equation*}
\Delta F(a,b)=\sum_{\varepsilon \in \{0,1\}^d}(-1)^{s(\varepsilon)}F(b+\varepsilon \ast (a-b)),
\end{equation*}

\Bin where for $\varepsilon=(\varepsilon_1,\cdots, \varepsilon_d)\in \{0,1\}^d$, $s(\varepsilon)=\varepsilon_1+\cdots+\varepsilon_d$

\bigskip \noindent An expanded version of that formula is :

\begin{equation*}
\Delta F(a,b)=\sum_{\varepsilon =(\varepsilon _{1},...,\varepsilon _{d})\in
	\{0,1\}^{d}}(-1)^{s(\varepsilon )}F(b_{1}+\varepsilon
_{1}(a_{1}-b_{1}),...,b_{d}+\varepsilon _{d}(a_{d}-b_{d})).
\end{equation*}

\noindent Let us try to understand the formula in a progressive way.\\

\bigskip \noindent \textbf{General rule of forming $\Delta F(a,b)$}. Let $a=(a_{1},...,a_{d})\leq b=(b_{1},...,b_{d})$ two points of $\mathbb{R}^{d}$ and let $F$ an arbitrary function
from $\mathbb{R}^{d}$ to $\mathbb{R}$. We form $\Delta F(a,b)$ in this way. First consider $F(b_{1},b_{2},...,b_{d})$ the value of $F$ at right endpoint $b=(b_{1},b_{2},...,b_{d})$ of the interval $]a,b].$ Next proceed to the replacement of each $b_{i}$ by $a_{i}$ by replacing exactly one of them, next two of them etc., and add the each value of $F$ at these points with a sign plus $(+)$ if the number of replacements is even and with a sign minus $(-)$ if the number of replacements is odd.\\

\noindent We recall the definition of distribution function on $\mathbb{R}$. 

\begin{definition} A function $F : \mathbb{R}^{d} \rightarrow \mathbb{R}$ is a distribution function \textit{(df)} on $\mathbb{R}^{d}$  if and only the two following conditions hold.\\
	
\noindent (a) $F$ assigns non-negative volumes to cuboids, that is $\Delta F(a,b)\geq 0$ for $a\leq b$.\\
	
\noindent (b) $F$ is right-continuous.\\
	
\noindent It is a probability distribution function \textit{pr.df} on $\mathbb{R}^{d}$\ if and only if the following three conditions are satisfied, where (c) is composed by two sub-conditions.\\
	
\noindent (a) $F$ assigns non-negative volumes to cuboids.
	
\noindent (b) $F$ is right-continuous.\\
	
\noindent (c) $F$ satisfies\\

 $(i)$
\begin{equation*}
\lim_{\exists i,1\leq i\leq d,t_{i}\rightarrow -\infty}F(t_{1},...,t_{d})=0
\end{equation*}
	
(ii)
\begin{equation*}
\lim_{\forall i,1\leq i\leq d,t_{i}\rightarrow +\infty
}F(t_{1},...,t_{d})=1.
\end{equation*}
\end{definition}

\bigskip\noindent The link between \textit{df}'s and Lebesgue-Stieltjes measures (LS-measures) is given by the following. We can associated to the \textit{df} $F$ a measure $\lambda_F$, called Lebesgue-Stieltjes measure associated to $F$, which is characterized by its values on the semi-algebra 

$$
\mathcal{S}=\{]a,b], \ a\leq b, \ (a,b) \in \overline{\mathbb{R}}^d\},
$$

\noindent which are

$$
\lambda_F(]a,b])=\Delta F(a,b).
$$

\Bin If $F$ is \textit{pr.df}, $\lambda_F$ is a probability measure. Conversely if $m$ is a measure on $\mathbb{R}^d$ such that

\begin{equation}
\forall x \in \mathbb{R}^d    \  F_m(x)= m(]-\infty, x])<\infty, \label{C}
\end{equation}

\Bin then $F_m$ is a \textit{df} (\textit{pr.df} if $m$ is a probability measure) such that $m=\lambda_{F_m}$.\\

\Bin \textbf{A2 - Spectrum and support}.\\

\Ni In this paper we need to introduce the notions of spectrum. First let us $\mathcal{O}$ as the class of all open sets in $\mathbb{R}^d$. We denote $\mathcal{N}(x)$the collection of neighborhoods of $x \in \mathbb{R}^d$. The spectrum of the \textit{df} $F$ is the following set 

$$
s(F)=\{x \in \mathbb{R}^d, \ \forall O \in \mathcal{N}(x), \ \lambda_F(O)>0\}.
$$

\Bin The point spectrum of $F$ is the set of atoms of  $\lambda_F$, that is

$$
ps(F)=\{x \in \mathbb{R}^d, \, \ \lambda_F(\{x\})>0\}.
$$

\Bin \noindent and the support of $F$ is the closure $\overline{ps(F)}$ of the point spectrum $ps(F)$.\\

\noindent \textbf{A3 - Moments}.  Let us define the class $\Gamma$ of a multi-indices in  $\mathbb{N}^d$, that is, all the row-vectors  $\alpha = (\alpha_1,\alpha_2,\cdots,\alpha_d)$ with  $\alpha_i \in \mathbb{N}$ for  $1\leq i\leq d$. Define the class of multi-index of level $\ell \in \mathbb{N}$.

$$
\Gamma(\ell)=\{\alpha=(\alpha_1,\alpha_2,\cdots,\alpha_d) \in \Gamma, \ |\alpha|\equiv \alpha_1+ \cdots + \alpha_d=\ell\}.
$$ 

\Bin For $u=(u_1,\cdots,u_d) \in \mathbb{R}^d$, we denote

$$
u^{\alpha}= \prod_{j=1}^{d} u_j^{\alpha_j}
$$

\noindent and the function $u \rightarrow u^\alpha$ is a polynomial of degree $|\alpha|$. Now we may define the moments of a \textit{df}.\\

\begin{definition} The moment of order $\alpha$ of a \textit{df} $F$ on $\mathbb{R}^d$ is the real number (whenever the integral exists) given by  

\begin{equation*}
\mu_\alpha = \int_{\mathbb{R}^d} \prod_{j=1}^{d} u_j^{\alpha_j} d\lambda_F(u_1,\cdots, u_d)\equiv \int_{\mathbb{R}^d} u^\alpha d\lambda_F(u).
\end{equation*}

\end{definition}

\Bin \noindent The moment problem we face in this paper amounts to the characterization of $F$ by all the field of moments $(\mu_\alpha)_{\alpha \in \Gamma}$, given they exist all.\\

\noindent \textbf{A4 - F-continuous interval}. First of all, a point $x=(x_1,\cdots,x_d)^t$ of $\mathbb{R}^d$ is a discontinuity point $x$, that is an element of the point spectrum $ps(F)$ of $F$ if and only if the boundary of $A_x=]-\infty, x]$ is not a $\lambda_F$-null set, i.e,

\begin{equation}
\lambda_F\left(\partial A_x \right)>0.
\end{equation}

\Bin We recall that   

$$
\partial A_x=\{y=(y_1,\cdots,y_d)^t \in \mathbb{R}^d, \forall j\in \{1,\cdots,d\} \ y_j\leq x_j, \ \exists j\in \{1,\cdots,d\} \ s.t. \ y_j=x_j\}.
$$

\Bin Further, for any an interval 

\begin{equation*}
(a,b)=\prod\limits_{i=1}^{d}(a_{i},b_{i})
\end{equation*}

\noindent of $\mathbb{R}^{d}$, we denote 

\begin{equation*}
E(a,b)=\{c=(c_{1},...,c_{d})\in \mathbb{R}^{d},\text{ }\forall 1\leq i\leq
d,(c_{i}=a_{i}\text{ ou }c_{i}=b_{i})\}.
\end{equation*}

\Bin By using the internal product ($8$) defined earlier, we have a compact form of $E(a,b)$ as 

\begin{equation} \label{sec_append_not}
E(a,b)=\{b+\varepsilon*(a-b), \varepsilon=(\varepsilon_1,...,\varepsilon_d) \in \{0,1\}^d\}.
\end{equation}

\bigskip \noindent  By definition, the interval $(a,b)$ is $F$-continuous if and only if $(a,b)$ is bounded and each element of $E(a,b)$ is a continuity point of $F$, that is

\begin{equation*}
\forall c\in E(a,b),\ \lambda_F(\partial ]-\infty ,c])=0.
\end{equation*}

\bigskip \noindent Let $\mathcal{U}(F)$ be the class of all $F$-continuous intervals. A key result which is very useful in weak convergence is the following proposition.\\

\begin{proposition} \label{cv.GFcontinuous}
Let $F$ be any probability distribution function on $\mathbb{R}^{d}$, $d\geq 1$. Then any open $G$ set in $\mathbb{R}^{d}$ is a countable union of $F$-continuous intervals of the form $]a,b]$ or $]a,b[$, where by definition, an interval $(a,b)$ is $F$-continuous if and only if, for any 

$$
\varepsilon=(\varepsilon_1, \varepsilon_2, ..., \varepsilon_d) \in \{0,1\}^d,
$$

\noindent the point 
$$
b+\varepsilon*(a-b)=(b_1+\varepsilon_1 (a_1-b_1), b_2+\varepsilon_2 (a_2-b_2), ..., b_d+\varepsilon_d (a_d-b_d))
$$ 

\noindent is a continuity point of $F$.
\end{proposition}

\Bin (see \cite{ips-wcia-ang}, Proposition 18, page 82 for a proof). A final consequence of that proposition is that any point $x=(x_1,\cdots,x_d)^t$ of $\mathbb{N}^d$ is limit of sequences of continuity points of $F$ from above and limit of sequences of continuity points of $F$ from below.\\

\noindent \textbf{B - An ordered version of Hahn-Banach theorem}.\\

\noindent Let us consider a linear space $E$ of real-valued functions $x$ defined on some space non-empty set $\Omega$ whose elements are denoted as

$$
x:\Omega\rightarrow \mathbb{R}.
$$ 

\bigskip\noindent Let $f$ be an element of the dual  space $E^{\prime}$ of $E$, that is, $f: E \rightarrow \mathbb{R}$ is a linear functional  (not necessary continuous). When we endow $E$ with the addition of functions and the external multiplication of functions by scalars and the following  partial order 

$$
\forall(x, y)\in E^2, (x\leq y) \Leftrightarrow \ (\forall t \in\Omega,  x(t)\leq y(t)),
$$

\Bin we can see that  $(E, +, ., \leq)$  is an $\mathbb{R}$-ordered linear space, that is, $(E, +, .)$  is an $\mathbb{R}$-linear space and the order relation is compatible with the linear structure, i.e.
       
$$
\forall(x, y, z)\in E^3, \    \   x\leq y \Leftrightarrow x+z \leq y+z
$$

\noindent and

$$
\forall (x, y)\in E^2, \forall \  \lambda \in\mathbb{R}_+\setminus \{0\}, \  (x \leq y)  \Leftrightarrow (\lambda x \leq \lambda y).
$$

\bigskip \noindent Given a non-empty subset $\Omega_0$ of $\Omega$ which may be equal to $\Omega$. We have the following definition.

\begin{definition} \label{partialNN} We say that $f \in E^{\prime}$ is $\Omega$-non-negative if and only if  

$$
(\forall x \in E \  and \ (\forall t\in \Omega_0, x(t) \geq 0) \Rightarrow (f(x) \geq 0),
$$

\Bin meaning that any function $x \in E$ which is non-negative on $\Omega_0$ has a non-negative image by $f$.
\end{definition}

\bigskip\noindent The following theorem is very similar to the Hahn-Banach theorem : given a linear sub-space $E_0$ of $E$ and given $f_0 \in E_0^{\prime}$ which is $\Omega$-non-negative, is it possible to extend $f \in E^{\prime}$ while preserving the $\Omega_0$-non-negativity? An affirmative response is given below.

\begin{theorem} \label{HBLIKE} Let $E$ be an ordered linear space of real-valued functions defined on some space non-empty set $\Omega$. Let $\Omega_0$ a non-empty subset of $\Omega$. Let $E_0$ be a sub-linear space of $E$. Let $f_0 \in E^{\prime}_0$ be $\Omega_0$-non-negative. Suppose that  $E_0$ has the following property:

\begin{equation}
\forall x \in E, \exists (x^{\prime}, x^{\prime\prime}) \in {E_0}^2,   \  \ x^{\prime} \leq x\leq x^{\prime\prime}  \  on \ \Omega_0, \label{ineq_00}
\end{equation}

\bigskip\noindent that is  

$$
\forall x\in E, \exists (x^{\prime}, x^{\prime\prime})\in E_0^2, \  (\forall t\in \Omega_0,  \ x^{\prime}(t) \leq x(t) \leq  x^{\prime\prime}(t)).
$$

\bigskip\noindent Then  $f_0$ is extensible to a linear functional on $E$ which is still $\Omega_0$-non-negative.
\end{theorem}  

\Ni \textbf{Proof}. We closely follow the proof of Hahn-Banach theorem which,by the way, is the approach used in \cite{shohat}. We notice that there is nothing to if  $E=E_0$ or for $E_0=\{0\}$,  $f_0=0$ and it is extended to $f=0$. So we proceed with $f_0 \neq 0 $ and $E\neq E_0 \neq \{0\}$. So there exists $x_0 \in E$ and $x_0 \notin E_0$. We consider the linear space spanned by $E_0$ and ${x_0}$ which is  

$$
E_1= E_0 + \mathbb{R} x_0 = \{ y= x+ \lambda x_0, x\in E_0, \ \lambda \in \mathbb{R} \}  
$$

\Bin We define on  $ E_1$ the functional 

$$
\forall y= x+\lambda x_0 \in E_1, f_1(y)=f_0(x)+\lambda a,
$$

\bigskip\noindent where $a$ is arbitrary real number and is taken as $f_1(x_0)$. For each $a$ fixed, $f_1$ is linear on $E_1$. $f_1$ is a extension of $f_0$ from $E_0$ to $E_1$, since any $y \in E_1$ is uniquely written as $y=x+\lambda x_0$ and then we have for $\lambda=0$,

$$
 f_1(y)=f_0(y) +0 a =f_0(x).
$$

\bigskip\noindent Now, the problem is how to choose $a=f_0(x_0)$ such that  $f_1$ is $\Omega$-non-negative. To do us, we begin by recalling the assumption

$$
A_1= \{ x^{\prime}\in E_0, x^{\prime}\leq x_0  \ on \   \Omega_0 \} \neq \emptyset \ and \  A_2 = \{ x^{\prime\prime} \in E_0, \ x^{\prime\prime}\geq x_0  \ on \   \Omega_0\} \neq \emptyset. 
$$

\Bin This implies that for any  $(x^{\prime}, x^{\prime\prime}) \in A_1 \times A_2$, $x^{\prime}\leq x_0 \leq x^{\prime\prime}$ on $\Omega_0$, and thus $(x^{\prime\prime}-x^{\prime})\geq 0$, on $\Omega_0$. Since $f_0$ is $\Omega_0$-non-negative, we have  $f_0(x^{\prime\prime}-x^{\prime})\geq 0$ [and hence 
$f_1(x^{\prime\prime}-x^{\prime})\geq 0$], that is

$$
\forall(x^{\prime}, x^{\prime\prime})\in A_1 \times A_2, \  f_0(x^{\prime})\leq f_0(x^{\prime\prime}) \ on \ \Omega_0.
$$

\bigskip\noindent Hence 

$$
\forall  x^{\prime}\in A_1,   \  f_0(x^{\prime}) \leq \inf_{ x^{\prime\prime}\in A_	2} f_0(x^{\prime\prime})  \ on \ \Omega_0,
$$

\bigskip\noindent Next, by taking the supremum on $x^{\prime}$, we have

$$
C_1=:\sup_{x^{\prime} \in A_1} f_0(x^{\prime})\leq \inf_{x^{\prime\prime} \in A_2} f_0(x^{\prime\prime})=:C_2.
$$

\bigskip\noindent Let us choose  $a \in [C_1, C_2]$. Let us show that for a such choice, $f_1$ will be $\Omega_0-$ non-negative.  Indeed, let  

$$
y=x+ \lambda x_0 \in E_1,
$$

\bigskip\noindent such that for $ t\in \Omega_0$, $y(t)=x(t)+ \lambda x_0(t) \geq 0$. We have to prove that $f_1(y) \geq 0$.  Let us discuss on the sign of $\lambda$.\\

\bigskip\noindent (a) Let $\lambda=0$. Here, for all $t \in \Omega_0$, $ y(t)= x(t) \geq 0$. So $f_0(x)=f_1(y) \geq 0$.\\

\noindent (b) Let $ \lambda>0$. Thus $(-x/\lambda) \leq x_{0}$  on $\Omega_0$. Thus $(-x/\lambda) \in A_1$. Hence

$$
f_1(x_0) \geq C_1= \sup_{x^{\prime} \in A_1}f(x^{\prime})\geq f_0(-\frac{x}{\lambda})
$$

\bigskip\noindent which leads to

\begin{eqnarray*}
f_1(x_0)- f_0\left(-\frac{x}{\lambda}\right)&=&{\frac{1}{\lambda}(f_0(x)+\lambda f_1(x_0))}\\  
	               &=&{\frac{1}{\lambda}f_1(y)\geq 0}
\end{eqnarray*}

\bigskip\noindent that is

$$
f_1(y) \geq 0.
$$

\Bin (c) Let $\lambda< 0$. Thus $(-x/\lambda) \geq x_0$. Thus $(-x/\lambda) \in A_2$. We use a similar argument to get 

$$
f_1(x_0)\leq C_2 = \inf_{x^{\prime\prime} \in A_2} f_0(x^{\prime\prime})  \leq f_0(-x/y)
$$ 

\bigskip\noindent and this leads to

$$
\frac{1}{\lambda}\left(f_0(x)+\lambda f_1(x_0)\right)={\frac{1}{\lambda}f(y)\leq 0}
$$  

\bigskip\noindent that is, since $\lambda<0$,

$$
f(y)\geq 0.
$$       

\bigskip\noindent We conclude that for $E_0 \subsetneq E$, we may extend $f_0$  to a bigger linear sub-space of at least on dimension, say $E_1$, as a linear and 
$\Omega_0$-non-negative functional.\\

\noindent For the second part, let us consider the class $\mathcal{A}$ of extensions of $f_0$ preserving $\Omega_0$-non-negativity. Let us denote them by $(f,A)$, meaning that $f: A\rightarrow \mathbb{R}$ is linear, $A$ subspace of $E$, $E_0\subsetneq A$ and $f_{|E_0}=f_0$ and $f$ is $\Omega_0$-non-negative. We say that $(f,A) \leq (f^{\prime},A^{\prime})$ if and only if

$$
(A\subseteq A^{\prime} \ \ and \ \ f^{\prime}_{|A}=f).
$$

\bigskip\noindent Clearly, this is an order relation. Let us exploit the first part. If $E_0\neq E$, there exists $x_0\neq E_0$ and $f_1 : E_1=E_0+\mathbb{R}x_0 \rightarrow \mathbb{R}$, $f_1\in \mathcal{A}$. If $E_1\neq E$, there exists $x_1 \in E\setminus E_1$ and we set $f_1 : E_2=E_1+\mathbb{R}x_1$, and we get $f_2\in \mathcal{A}$.\\

\noindent \textbf{Either}, we stop at some $n$ with $E_n=E$, and the proof is finished \textbf{or} we continue infinitely. But, by construction, we have 

$$
(f_0,E_0)\leq (f_1,E_1)\leq (f_2,E_2)\leq \cdots \leq (f_j,E_j) \cdots
$$

\bigskip\noindent So the class $\lbrace (f_j,E_j), j\geq 0 \rbrace$ is a chain. The  Zorn's lemma says that it has a maximal element. It is not difficult to see that this maximal element is $(f_{\infty},E_{\infty})$ with 

\begin{equation*}
 \left\{
	\begin{array}{lll}
 E_{\infty}=\bigcup_{j\geq 0}E_j\\
 \\
\forall x \in E_{\infty},\ f_{\infty}(x)=f_j(x), \ \text{for} \ x \in E_j\\
    \end{array}.
\right.
\end{equation*}

\bigskip\noindent Since the $(E_j)_{j\geq 0}$ is an increasing sequence (w.r.t to the inclusion), 

$$
E_{\infty}=\bigcup_{j\geq 0}E_j
$$ 

\Bin is a linear sub-space of $E$. Let us see that the definition is coherent. Indeed, let us suppose that $x \in E_{\infty}$ belongs two distinct spaces $E_{j_1}$ and $E_{j_2}$, $j_1\geq 0$ and $j_2\geq 0$. Without loss of generality, we can suppose that $j_1<j_2$.  Hence, we have 

$$
(f_{j_1},E_{j_1})\leq (f_{j_2},E_{j_2}).
$$

\bigskip\noindent and thus,

$$
f_{j_1}(x)=f_{{j_2}| E_{j_1}}(x)=f_{j_2}(x).
$$

\bigskip\noindent We may take $f_{\infty}(x)$ as $f_{j}(x)$ for any $j\geq 1$ such that $x\in E_j$. All these values are equal by the previous formula. So, the definition of $f_{\infty}$ is coherent.\\

\noindent The mapping $f_{\infty}$ is linear since for $x\in E_{\infty}$, $y\in E_{\infty}$, there exist $j_1$ and $j_2$ (say $j_1\leq j_2$) such that $x\in E_{j_1}$ and 
$y \in E_{j_2}$. So $(x,y) \in E_{j_2}$. For $(\alpha,\beta)\in \mathbb{R}^2$, $\alpha x +\beta y \in E_{j_2}$

\begin{eqnarray*}
f_{\infty}(\alpha x +\beta y)&=&f_{j_2}(\alpha x +\beta y)=\alpha f_{j_2}( x )+\beta f_{j_2}(y)\\
&=&\alpha f_{\infty}(x)+ \beta f_{\infty}(y).
\end{eqnarray*}

\bigskip\noindent We have $E_0 \subseteq E_{\infty}$ obviously and for all $j\geq 0$, for all $x\in E_j$

$$
f_{\infty}(x)=f_{j}(x).
$$ 

\bigskip\noindent So, $f_{\infty|E_j}(x)=f_j(x)$ and hence $f_{\infty|E_0}(x)=f_0(x)$. We also have that $f_{\infty}$ is $\Omega_0$-non-negative. Indeed for $x \in E_{\infty}$, $x \geq 0$ on $\Omega_0$, we have for $x \in E_j$, $f_{\infty}(x)=f_j(x)\geq 0$.\\

\noindent So $f_{\infty} $ belongs to $\mathcal{A}$ and dominates all elements of $\mathcal{A}$. Hence

$$
(f_{\infty}, E_{\infty})=\max \lbrace (f_{j}, E_{j}), j\geq 0 \rbrace.
$$

\bigskip\noindent We necessarily have $E_{\infty}=E$. Indeed if $E_{\infty}\subsetneq E$, we might use the first part and set 
$0\neq x_{\infty} \in E \setminus E_{\infty}$ and we obtain a greater extension $f_{\infty}^\ast$ preserving the $\Omega_0$-non-negativity, defined on 
$E_{\infty}^\ast=E_{\infty}+\mathbb{R}x_{\infty}$, which is impossible. $\blacksquare$\\
            
\newpage \section{The moment problem in Probability Theory of $\mathbb{R}$} \label{mp_r_proba} 

\noindent Suppose that we have a probability measure $\rho$ on $\mathbb{R}$ having moments of all orders $(m_n)_{n\geq 1}$, with $m_0=1$, as in Formula \eqref{mp_01}. The question is whether the sequence characterizes the measure $\rho$ in the following form : If $(m_n)_{n\geq 0}$, with $m_0=1$, are the moments of two measures 
$\rho_1$ and $\rho_2$ on $\mathbb{R}$, do we have $\rho_1=\rho_2$? We have the particular answer as follows.\\
 
\textbf{(I) - A sufficient condition for the moments to determine the probability measure}.\\

\begin{theorem} \label{mp_proba_r} Let $\rho$ be a probability measure on $\mathbb{R}$ having moments of all orders $(m_n)_{n\geq 1}$, with $m_0=1$. Suppose that the Cauchy radius exists and is not zero, i.e.,

$$
R=\lim_{n \rightarrow +\infty} |n!/m_n|^{1/n}>0,
$$

\Bin or the series $\sum_{n=0}^{+\infty} m_n x^n/n!$ has a positive radius of convergence.\\
 
\Bin Then the moments determine $\rho$.
\end{theorem}

\Bin The simple tool of Cauchy's rule for convergence of functional series is used here. Let us just make a recall. Let us consider a sequence of real numbers $(a_n)_{n\geq 0}$. Suppose that $|1/a_n|^{1/n}\rightarrow r>0$. Then for $|x|<r$, such that
$0<\varepsilon=1-|x/r|>0$. We have

\begin{eqnarray*}
|a_n x^n|&=&\biggr(|x| |a_n|^{1/n}\biggr)^n\\
&=&\biggr(\left|\frac{x}{r} \right| \biggr[\left|r |a_n|^{1/n}\right|\biggr]\biggr)^n.
\end{eqnarray*}

\Bin Since the term between the brackets converges to one, it is less that $(1-\varepsilon/2)^{-1}>1$ for $n$ large enough, say $n\geq n_0$. We get

$$
|a_n x^n|\leq \left(\frac{1-\varepsilon}{1-\varepsilon/2}\right)^n.
$$

\Bin since $0<(1-\varepsilon)/(1-\varepsilon/2)<1$, the series $\sum_{n} a_n x^n$ converges for all $|x|<r$. Similarly, we prove that the series $\sum_{n} a_n x^n$
diverges for $|x|>r$. We are going to use that rule below based on arguments in \cite{billinsgleymp}, page 388.\\

\noindent \textbf{Proof of Theorem \ref{mp_proba_r}}. Let us denote by $\psi$ the characteristic function of $\rho$. The Taylor-Lagrange formula (see \cite{valiron}, p. ??) for  the complex exponential function gives : for $(x,t,h)\in \mathbb{R}^3$, $n\geq 1$,

$$
e^{ihx}=\sum_{j=0}^{n}\frac{(ihx)^j}{j!} +\frac{(i xh)^{n+1} e^{(i\theta xh)}}{(n+1)!}, \ |\theta|<1.
$$

\Bin This leads to (since $e^{itx}$ has norm one)

$$
\left|e^{itx} \left(e^{ihx}-\sum_{j=0}^{n}\frac{(ihx)^j}{j!}\right)\right|\leq \frac{|xh|^{n+1}}{(n+1)!},
$$

\Bin which yields

$$
\left|e^{i(t+h)x} -\sum_{j=0}^{n}\frac{h^j}{j!} (ix)^j e^{itx} \right|\leq \frac{|h|^{n+1}}{(n+1)!} |x|^{n+1}.
$$

\Bin By integrating the three members of that double inequality with respect to $\rho$ and by identifying $\int (ix)^j e^{itx} \rho(x)$ as the derivative of $\psi$ at $j$, we get

\begin{equation}
\left|\psi(t+h)- \sum_{j=0}^{n}\frac{h^j}{j!} \psi^{(j)}(t) \right| \leq \frac{|h|^{n+1}}{(n+1)!} \mu_{n+1}, \label{mp_tool_01}
\end{equation}

\Bin where $\mu_{j}$ is the absolute moment of order $n\neq 1$ given by

$$
\forall j\geq 1, \ \mu_{j}=\int \ |u|^j \ d\rho(u).
$$

\Bin Now, under the hypotheses, we can find $r$ and $s$ such that $0<r<s<1$ and $\sum_{j=0}^{+\infty} m_j s^{j}/j!$ converge. Hence by the properties of convergent series, $m_j s^{j}/j!\rightarrow 0$ and $m_j r^{j}/j!\rightarrow 0$ as $j\rightarrow +\infty$. Further, $2 e^{\log j + (2j-1) \log (r/s)} \rightarrow -\infty$ (since $0<r/s<1)$ and then $2 e^{\log j + (2j-1) \log (r/s)} <s$ for $j$ large enough, say $j\geq j_0$, which is

$$
2j r^{2j-1} <s^{2j}, \ for \ j\geq j_0,
$$

\Bin which, combined with the inequality $|a|^{r_1} \leq 1 +|a|^{r_2}$ valid for $0<r_1\leq r_2$, leads to

\begin{eqnarray*}
\frac{|x|^{2j-1}r^{2j-1}}{(2j-1)}&\leq& \frac{r^{2j-1}}{(2j-1)!}+\frac{|x|^{2j}r^{2j-1}}{(2j-1)!}\\
&\leq &\frac{r^{2j-1}}{(2j-1)!}+\frac{|x|^{2j}s^{2j-1}}{(2j-1)!(2j)}\\
&\leq &\frac{r^{2j-1}}{(2j-1)!}+\frac{|x|^{2j}s^{2j}}{(2j)!}.
\end{eqnarray*}

\Bin By integration with respect to $\rho$, we get

$$
\frac{\mu_{2j-1}r^{2j-1}}{(2j-1)!}\leq \frac{r^{2j-1}}{(2j-1)}+\frac{\mu^{2j}s^{2j-1}}{(2j)!}.
$$

\Bin So $\mu_{2j-1}r^{2j-1}/(2j-1)! \rightarrow 0$. We already have $\mu_{2j}r^{2j}/(2j)!=m_{2j}r^{2j}/(2j)! \rightarrow 0$. So, the convergence covers odd and even terms. We arrive at 

$$
\mu_{n+1}r^{n+1}/(n+1)! \rightarrow 0 \ as \ n\rightarrow 0.
$$

\Bin We apply this to the bound in Formula \ref{mp_tool_01} to get

\begin{equation}
\forall t \in \mathbb{R}, \ \forall |h|\leq r, \ \psi(t+h)= \sum_{j=0}^{+\infty}\frac{h^j}{j!} \psi^{(j)}(t). \label{mp_tool_02}
\end{equation}

\Bin We conclude as follows. Let us suppose that another probability measure has the same moments $(m_n)_{n\geq 1}$ with characteristic function $\psi_1$. By taking For $t=0$, we get that $\psi$ and $\psi_1$ coincide on $[-r,r]$. Let us  show we may extend that equality to all interval $[sr, \ (s+1)r]$, $s\geq 1$. We begin by preceeding for $s=1$. We say that $\psi$ and $\psi_1$ have the same derivative functions on $]0,r[$ and $\psi^{(j)}(r/2)=\psi^{(j)}_1(r/2)$ for all $j\geq 1$ in particular. By taking $t=r/2$, Formula \eqref{mp_tool_02} shows that $\psi$ and $\psi_1$ are equal on $[r/2, 3r/2]$ and hence $\psi^{(j)}(r)=\psi^{(j)}_1(r)$ for all $j\geq 1$. Now using Formula \eqref{mp_tool_02} extends the equality on $[r, 2r]$. By proceeding so forth and by handling intervals $[-(s+1)r, \ -sr]$ in the same way, we get the desired equality on $\mathbb{R}$ by induction. $\blacksquare$\\

\newpage

\noindent \textbf{(II) - Application to weak convergence}.\\

\noindent We get the following criteria of convergence.\\

\begin{theorem} \label{wc_moments} Let $X_n : (\Omega, \mathcal{A}_n, \mathbb{P}_n) \rightarrow \mathbb{R}$, $n\geq 1$, be a sequence of random variables and $X_{\infty} : (\Omega_{\infty}, \mathcal{A}_{\infty}, \mathbb{P}_{\infty}) \rightarrow \mathbb{R}$ be another random variable. Let us suppose that the $X_n$'s and $X_{\infty}$ have moments of all orders and that the probability law of $X_{\infty}$ is determined by its moments and

$$
\forall j\geq 1, \ \mathbb{E}_{\mathbb{P}_n} X_n^j \rightarrow \mathbb{E}_{\mathbb{P}_{\infty}} X_{\infty}^j \ as \ n\rightarrow +\infty.
$$

\Bin Then $X_n$ weakly converges to $X_{\infty}$ as $n\rightarrow +\infty$, i.e., $X_{n}\rightsquigarrow X_{\infty}$.
\end{theorem}

\Bin \textbf{Proof}. Since the sequence $\mathbb{E}_{\mathbb{P}_n} X_n^2$ converges, it is bounded, say by $C$. For any $\varepsilon>0$, for $k>0$ and $C/k2<\varepsilon$, we apply the Markov inequality to get

$$
\mathbb{P}_n(|X_n|\geq k)=\mathbb{P}_n(X_n^2\geq k^2)\leq C/k^2<\varepsilon,
$$  

\Bin that is, there exists a compactum $K=[-k,k]$ of $\mathbb{R}$ such that

$$
\liminf_{n\rightarrow +\infty} \mathbb{P}_n(X_n \in K)>1-\varepsilon.
$$

\Bin So the sequence $(X_n)_{n\geq 1}$ is asymptotically tight and by Prohorov's theorem, every sub-sequence of $(X_n)_{n\geq 1}$ contains a weakly convergent sub-sequence (see Theorem Prohorov-Helly Bray in \cite{ips-wcia-ang}, Section 3, Sub-section 3). Now let $f : \mathbb{R} \rightarrow \mathbb{R}$ continuous and bounded. The sequence $s_n(f)=\mathbb{E}_{\mathbb{P}_n}f(X_n)$ is bounded (by the bound of $f$). So, it contains converging sub-sequence $s_{n_k}(f)$ to $s(f)$. But, by Prohorov's theorem, $X_{n_k}$ contains a sub-sequence $X_{n_k(\ell)}$ weakly converging, say to $Z$ of probability measure $\mu$. So 

$$
s(f)=\int f \ \ d\mu.
$$

\Bin Let us use the Skorohod theorem (see \cite{wichura} ) to have $X^{\ast}_{n_k(\ell)}=_d X_{n_k(\ell)}$ and $Z^{\ast}=_d Z$ on the same probability space with 
$X^{\ast}_{n_k(\ell)}$ converges \textit{a.s.} to $Z^{\ast}$. For any $r\geq 1$ fixed, $\mathbb{E} (X^{\ast}_{n_k(\ell)})^{4r}$ is bounded and hence, for any $r\geq 1$, $(X^{\ast}_{n_k(\ell)})^{r}$ is uniformly and continuously integrable and converges to $(Z^{\ast})^r$. By Theorem 16.4 in \cite{billinsgleymp}, page 218 , $(Z^{\ast})^r$ is integrable and $\mathbb{E} (X^{\ast}_{n_k(\ell)})^{r}$ converges to $\mathbb{E}(Z^{\ast})^r$. By getting back to our original random variables, we get

$$
\mathbb{E}X_{n_k(\ell)}^r \rightarrow \mathbb{E} Z^r.
$$

\Bin Since $\mathbb{E}X_{n_k(\ell)}^r$ converges to $\mathbb{E} X^r_{\infty}$, we get that $X_{\infty}$ and $Z$ have the same moments (which determine the probability law of $X_{\infty}$), we conclude that $\rho=\mathbb{P}_Z=\mathbb{P}_{X_\infty}$. Hence

$$
s(f)=\int f \ \ d\mathbb{P}_{X_{\infty}}.
$$

\Bin We conclude that any sub-sequence of $s_n(f)$  contains a sub-sequence converging to $s(f)=\int f \ d\mathbb{P}_{\infty}$ for any bounded and continuous function $f$. Thus, $X_n \rightsquigarrow X_{\infty}$. $\blacksquare$\\

\newpage \section{Solution of the moment problem on $\mathbb{R}$ and application to weak convergence} \label{mp_r_shohat}

\noindent Here, we use simpler notations. Let be given the sequences $(m_n)_{n\geq 0}$ with $m_0=1$. Let $\mathcal{P}$ the linear space of all polynomials. A non-zero polynomial $P$ is associated with coefficients $(x_n)_{n\geq 0}$, where all the $x_n$'s vanish beyond some integer $d$ for which $x_d\neq 0$, the number $d$ being its degree. For sake of simplicity, we use the representation $P\equiv (x_n)_{n\geq 0}$ and use infinite sums with in mind the fact that only a finite number of the sum is non-zeros :

$$
\forall u \in \mathbb{R}, \ P(u)=\sum_{n\geq 0} x_n u^n.
$$ 

\noindent We define the linear functional $\mu$ as follows :

$$
\forall P\equiv (x_n)_{n\geq 0} \in \mathcal{P}, \ \mu(P)=\sum_{n\geq 0} x_n m_n.
$$

\Bin That functional is well-defined and is linear. Here is the solution of the moment problem on $\mathbb{R}$.

\begin{theorem} \label{shohat_R} Given a non-empty closed subset $S_0$ of $\mathbb{R}$, there exists a probability measure $\rho$ associated to a \textit{df} $F$ such that : (a) $supp(F)\subset S_0$ and (b) for all $n\geq 0$,

$$
m_n= \int u^n \ d\rho(u)
$$

\Bin if and only if : (c) $\mu$ is $S_0$-non-negative, i.e., if $\mathcal{P} \ni P$ satisfies : $P(u)\geq 0$ for all $u\in S_0$, then $\mu(P)\geq 0$.
\end{theorem}

\Bin \textbf{Proof}. We are going to provide a detailed proof.\\

\noindent \textbf{Let us begin by proving that (a) and (b) imply (c)}. For any polynomial $P=(x_n)_{n\geq 0}$ $S_0$-non-negative, we have

$$
\mu(P)=\sum_{n\geq 0 1} x_n \left(\int u^n \ dF(u)\right)=\int \left(\sum_{n\neq 1} x_n \ u^n \right) dF(u)=\int P(u) dF(u),
$$

\noindent where we were able to interchange summation and integration symbols since only a finite number of terms of the summation are non-zero. But, we have

$$
\mu(P)= \int P(u) dF(u) = \int_{S_0^c} P(u) dF(u) + \int_{S_0} P(u) dF(u).
$$

\Bin But, on $\mathbb{R}$, the support $supp(F)$ and the spectrum $s(F)$ coincide and since $S_0^c \subset supp(F)^c$, we have 

$$
\int_{S_0^c} P(u) dF(x)=0
$$

\Bin and we get 

$$
\mu(P)= \int P(u) dF(u) = \int_{S_0} P(u) dF(u).
$$

\Bin which is non-negative whenever $P$ is $S_0$-non-negative.\\

\noindent \textbf{Let us prove that (c) implies (a) and (b)}. Let us proceed with three steps.\\

\noindent \textit{Step 1. Construction of $\rho$}. \label{step1} Let us consider the class $E$ of functions $f : \mathbb{R} \rightarrow \mathbb{R}$ bounded of linear combinations of functions of the form $A u^{2r} + B$, where $A\geq 0$, $B\geq 0$, $r \in \mathbb{N}$. In other words $f \in E$ if and only if it is bounded by a function of the form

\begin{equation}
g=\sum_{i=1}^{p} A_i u^{2r_i} + B_i, \ p\geq 1, \ (A_i,B_i,r_i) \in \mathbb{R}_+ \times \mathbb{R}_+ \times \mathbb{N}. \label{eq_01}
\end{equation}

\noindent We set $E_0=E \cap \mathcal{P}$ as the subclass of $E$ restricted to polynomials. It is clear that for a function $g$ as in Formula \eqref{eq_01}, $-g$ and $g$ belong to $E_0$ and hence : 

\begin{equation}
\forall f \in E, \ \exists (f_1, f_2) \in E_0^2, \  \ f_1 \leq f \leq f_2 \ on \ \mathbb{R}. \label{ineq_01}
\end{equation}

\Bin We may apply Theorem \ref{HBLIKE} since $E_0$ is a sub-linear space of $E$, $\mu$ is an $S_0$-non-negative linear functional defined on $E_0$ and Condition \eqref{ineq_00} of Theorem is true through Formula \eqref{ineq_01}. So $\mu$ est extensible on $E$ to an $S_0$-non-negative linear functional, still denoted by $\mu$. For any subset $C$ of $\mathbb{R}$, $f=1_C$ is bounded by $g=1= 0 \times u^2 + 1$ so that $1_C \in E$. So define the mapping $m$ on the class $\mathcal{B}(\mathbb{R})$ of Borel sets of $\mathbb{R}$ by

$$
\forall C \in \mathcal{B}(\mathbb{R}), \ m(C)=\mu(1_C).
$$

\Bin The mapping is clearly additive. For any $C \in \mathcal{B}(\mathbb{R})$, $1_C\geq 0$ on $\mathbb{R}$ and thus on  $S_0$, we have by $S_0$-non-negativity of $\mu$, that $\mu(C)=\mu(1_C)\geq 0$. As well, for $(C_1,C_2) \in \mathcal{B}(\mathbb{R})^2$,  $C_1 \subset C_2$ implies  $1_{C_1} \leq 1_{C_2}$ on $\mathbb{R}$ and hence on $S_0$ and by $S_0$ non-negativity of $\mu$, $m(C_1)\leq m(C_2)$. Finally

$$
\forall C \in \mathcal{B}(\mathbb{R}), \ m(C) \leq \mu(1_{\mathbb{R}}=\mu(1)=m_0=1.
$$

\noindent We conclude that $m$ is a finite and non-negative additive mapping on $\mathcal{B}(\mathbb{R})$. The mapping should be a measure if we could prove that it is $\sigma$-sub-additive or continuous at $\emptyset$, that is $m(A_n)\downarrow 0$ if $A_n\downarrow \emptyset$ as $n\downarrow +\infty$. But it seems very difficult to prove that. So we are going to use the same method as in \cite{shohat} but in the modern frame of Measure Theory.\\

\noindent We define the function $F_0(x)=m(]-\infty, \ x])$, $x \in \mathbb{R}$. We are not sure that it is right-continuous. So we work with

$$
F(x)=\lim_{h \searrow 0} F_{0}(x+h), \ x \in \mathbb{R}.
$$

\noindent The limits exist by the monotonicity of $F_0$ and the function $F$ is right-continuous and assigns to intervals $]a,b]$ non-negative lengths, that is 
$\Delta F(a,b)=F(b)-F(a)\geq 0$. Hence $F$ is a distribution function. Let us denote $\rho=\lambda_F$ the Lebesgue-Stieltjes measure associated with $F$.\\

\noindent It is useful to remark that $F_0$, as a monotone function, has at most a countable number of discontinuity, so that 

$$
\int 1_{]a,b]} \ d\rho =F(b)-F(a)=F_0(b)-F_0(a),
$$ 

\Bin except, eventually, for at most a countable number of pairs $(a,b)$. Now let us check that : Any non-negative and increasing or decreasing function $f \in E$ is $\rho$-integrable and we have

\begin{equation}
0\leq \int f d\rho \leq \mu(f). \label{boundedness}
\end{equation}

\Bin  Let us finish this step by proving the above claims. We suppose that $f$ is increasing. By definition of the integral with respect to $\rho$, the integral of $f$ is the monotone limit of integrals of a sequence $(g_n)_{n\geq 1}$ of elementary functions, each of them having the following form

$$
g= \sum_{j=1}^{p} \alpha_j 1_{(a_j\leq f <b_j)}, \ p>1, \ p \ finite, \ (\alpha_j)_{1\leq j \leq p} \subset \mathbb{R}_{+}
$$

\Bin with $0\leq g \leq f$. But we have

\begin{eqnarray*}
\int g \ d\rho &=& \sum_{j=1}^{p} \alpha_j \rho(a_j\leq f <b_j)\\
&=& \sum_{j=1}^{p} \alpha_j \rho([f^{-1}(a_j), f^{-1}(b_j)[)\\
&=& \sum_{j=1}^{p} \alpha_j \{F(f^{-1}(b_j)+0)-F(f^{-1}(a_j)-0)\},
\end{eqnarray*}

\noindent where $F(x+0)$ and $F(x-0)$ are the left and the right limit of $F$ at $x$ respectively. The boundaries $a_j$ and $b_j$ can be chosen as continuity points of $F_0$ (which still are continuity points of $F$), the only requirement being that the modulii $b_j-a_j$ be small enough. Hence

\begin{eqnarray}
\int g \ d\rho &=&\sum_{j=1}^{p} \alpha_j \{F_0(f^{-1}(b_j))-F_0(f^{-1}(a_j)\} \notag\\
&=&\sum_{j=1}^{p} \alpha_j m(a_j\leq f <b_j)= \mu \left(\sum_{j=1}{p} \alpha_j 1_{(a_j\leq f <b_j)}\right) \notag\\
&=& \mu(g)\leq \mu(f). \notag
\end{eqnarray}

\Bin So, for all $n\geq 1$, 

$$
0\leq \mu(g_n) = \int g_n \ d\rho.
$$

\noindent At the limit, we have $\int f \ d\rho \leq \mu(f)$. Hence $f$ is $\rho$-integrable and its integral is bounded by $\mu(f)$. The proof is easily adapted for $f$ decreasing. Let us give an example. For each function $\ell_n(u)=u^n$, $\ell_n^{+}$ and $\ell_n^{-}$ are still in $E$ and the bound given above applies to them and we finally have 

$$
\left|\int \ell_n(u) \ d\rho(u)\right| \leq \int \ell_n(u)^{+} \ d\rho(u)+ \int \ell_n(u)^{-} \ d\rho(u)\leq \mu(\ell_n^{+})+\mu(\ell_n^{-})=\mu(|\ell_n|),
$$

\Bin we have the following

\begin{fact} \label{factIntegrability} For any $n\geq 0$, the function $\ell_n(u)=u^n$ of $u\in \mathbb{R}$ is $\rho$-integrable and 

\begin{equation}
\left| \int \ell_n d\rho \right| \leq \mu(|\ell_n|). \label{boundednessMn}
\end{equation}

\end{fact}

\Bin  \textit{Step 2. $s(F) \subset S_0$}. Let us prove that $S_0^c \subset s(F)^c$. Let $x \in S_0^c$, which is an open set. So there exists an interval $]a,b[$ such that $x \in ]a,b[$ and $]a,b]\subset S_0^c$. The number $a$ and $b$ can be taken as continuity points of $F_0$. Since $1_{]a,b]}=0$ on $S_0$, i.e., $1_{]a,b]}$ is non-positive on $S_0$, we have $\mu(1_{]a,b]})\leq 0$ and since $\mu(1_{]a,b]})\geq 0$, we have 

$$
0=\mu(1_{]a,b]})=F_0(b)-F_0(a)=F(b)-F(a)=\rho(]a,b])\geq \rho(]a,b[).
$$

\Bin Since $G=]a,b[$ is an open neighborhood of $x$ such that $\rho(G)=0$, we conclude that $x \notin s(F)$. Let us move to the last step.\\

\noindent \textit{Step 3. $\rho$ has the desired moments}.\\

\noindent Let $n\geq 1$. Let us show that

$$
m_n=\int u^n \ d\rho(u).
$$

\Bin Let $\varepsilon \in ]0,  1[$ be fixed. Let $K$ be a positive integer such that $1/K \leq \varepsilon$ (and thus $K\geq 1$). Hence for an positive integer $r$ such that $2r-n-1\geq 1$, we have for $u \notin ]-K,K[$,

\begin{eqnarray*}
|u|^n&=& u^{2r} \frac{1}{|u|^{2r-n}}\\
&\leq& u^{2r} \frac{1}{K^{2r-n}}=\frac{u^{2r}}{K} \frac{1}{K^{2r-n-1}}\\
&\leq& \varepsilon u^{2r}.
\end{eqnarray*}

\Bin We conclude that for $K$ such that $1/K \leq \varepsilon$, for $u \notin ]-K,K[$

\begin{equation}
|u|^n \leq \varepsilon u^{2r} \leq  u^{2r}.\label{boundE}
\end{equation}

\Bin Now, the function $\ell_n(u)=u^n$ of $u\in \mathbb{R}$ is uniformly continuous on $I_K=]-K, \ K]$. Let us fix $\eta>0$ and let us divide $]-K,\ K]$ into a finite number $p$ of intervals $]a_h,b_h]$ such that the variation of $\ell_n$ over $]a_h,b_h]$ is less than $\eta$. It is possible to choose the $a_h$'s and the $b_h$'s as continuity points of $F_0$ (and hence of $F$). [To do that, we may divide each intervals into two at the middle and to move each $a_h$ and $b_h$ very slightly to be a continuity point. The variation of $\ell_n$ over the new intervals remain is less than $\eta$].\\

\Ni Let us define an elementary function $\ell_{p,n}$ by choosing $u_{(h)}$ from each interval $]a_h,b_h]$ as follows

$$
\ell_{p,n}(u)=\sum_{j=1}^{p} \ell_n(u_{(h)}) 1_{]a_h,b_h]}, \ u \in \mathbb{R}.
$$

\Bin We have

\begin{eqnarray}
\mu(\ell_{p,n})&=&\sum_{j=1}^{p} \ell_n(u_{(h)}) \mu\left(1_{]a_h,b_h]}\right) \label{ellRho}\\
&=&\sum_{j=1}^{p} \ell_n(u_{(h)}) (F_0(b_h)-F_0(a_h) \notag\\
&=&\sum_{j=1}^{p} \ell_n(u_{(h)}) (F(b_h)-F(a_h)) \notag\\
&=&\int \ell_{p,n} \ d\rho.  
\end{eqnarray}

\Bin we notice that $\ell_{p,n}$ is null on $]-K, \ K]^c$. By using Formula \eqref{boundE} and the continuity modulus of $\ell_n$ over $]-K, \ K]$, we have

\begin{eqnarray*}
|\ell_n(u)-\ell_{p,n}(u)| &\leq& |\ell_n(u)-\ell_{p,n}(u)| 1_{]-K, \ K]}+ |\ell_n(u)-\ell_{p,n}(u)| 1_{]-K, \ K]^c}\\
&\leq& \eta 1_{\mathbb{R}}+ \varepsilon u^{2r},
\end{eqnarray*}

\Bin i.e., for all $u\in \mathbb{R}$,

\begin{equation}
\ell_{p,n}(u) - \eta - \varepsilon u^{2r} \leq \ell_n(u) \leq +\ell_{p,n}(u) + \eta + \varepsilon u^{2r}.
\end{equation}

\Bin By applying $\mu$ to that ordering on $\mathbb{R}$ (and hence on $S_0$) and by using Line \eqref{ellRho}, we get

\begin{equation}
\int \ell_{p,n} \ d\rho - \eta  \mu(1_{\mathbb{R}}) - \varepsilon m_{2r} \leq m_n \leq \int \ell_{p,n} \ d\rho  +\eta \mu(1_{\mathbb{R}}) +\varepsilon m_{2r}.
\end{equation}

\Bin We notice that

$$
\int \ell_{p,n} \ d\rho=\int 1_{]-K, \ K]} \ell_{p,n} \ d\rho.
$$ 

\Bin On $]-K, \ K]$, $\ell_{p,n} \rightarrow \ell_n$ and bounded by $|\ell_n|$ which is integrable by Fact \ref{factIntegrability}. By letting $\eta \downarrow 0$, we will have $p\rightarrow +\infty$ and the dominated convergence theorem, as 

$$
\int \ell_{p,n} \ d\rho \rightarrow \int 1_{]-K, \ K]} \ell_{n} \ d\rho.
$$

\Bin and hence

\begin{equation}
\int 1_{]-K, \ K]} \ell_{n} \ d\rho - \varepsilon m_{2r} \leq m_n \leq \int 1_{]-K, \ K]} \ell_{n} \ d\rho + \varepsilon m_{2r}.
\end{equation}

\Bin For $\varepsilon>0$ fixed, we can let $K\uparrow +\infty$, $1_{]-K, \ K]} \ell_{n} \rightarrow \ell_n$ while being dominated by the integrable function $|\ell_n|$ and hence

\begin{equation}
\int \ell_{n} \ d\rho - \varepsilon m_{2r} \leq m_n \leq \int \ell_{n} \ d\rho + \varepsilon m_{2r}.
\end{equation}

\Bin Now, we may let $\varepsilon \rightarrow 0$ to get

$$
m_n=\int \ell_{n} \ d\rho. \blacksquare
$$





\begin{thebibliography}{99}
\bibitem[lo\`eve (1997)]{loeve} Lo\`{e}ve, Michel.(1997). \textit{Probability Theory I}. Springer-Verlag, 4th Edition.	
\bibitem[Lo(2017b)]{ips-mestuto-ang} Lo, G. S. (2017) Measure Theory and Integration By and For the Learner.
 \text {SPAS Books Series}. Saint-Louis, Senegal - Calgary, Canada. Doi :  http://dx.doi.org/10.16929/sbs/2016.0005,
ISBN : 978-2-9559183-5-7
\bibitem[Lo (2016)]{ips-probelem-ang} Lo, G.S.(2016). \textit{A Course on Elementary Probability Theory}. SPAS Editions. Saint-Louis, Calgary, Abuja. Doi : 10.16929/sbs/2016.0003.
\bibitem[Lo (2018)]{ips-mfpt-ang} Lo, G.S.(2018). \textit{Mathematical Foundation to Probability Theory}. Spas Textbooks Series.
\bibitem[Lo (2017)]{ips-wcia-ang} Lo, G.S.(2018). \textit{Weak Convergence (IA) - Sequences of Random Variables}. SPAS Book Series, Calgary, Alberta, Saint-Louis (Sénégal).
\bibitem[Shohat and Tamarkin (1943)]{shohat} Shohat J.A. and Tamarkin J.D (1943) The problem of moment, Mathematical Surveys and monographs. Volume I. American Society of Mathematics (Re-edited in 1950, 1963 and 1970).
\bibitem[Billingsley(1995)]{billinsgleymp}  Patrick Billingsley (1995). \textit{Probability and Measure}. Wiley. Third Edition.
\bibitem[Billingsley (1968)]{billingsley} Billingsley, P.(1968). \textit{Convergence of Probability measures}. John Wiley, New-York.
\bibitem[Gutt (2005)]{gutt}  (2005). Allan Gutt (2005). \textit{Probability Theory : a graduate course}. Springer.
\bibitem[Valiron (1941)]{valiron}  (1941). Valiron George (1941). \textit{Th\'eorie des fonctions}. Masom, Paris.
\bibitem[Wichura (1996)]{wichura}  (1996). 
\end{thebibliography}
\end{document}